\documentclass{amsart}

\usepackage{amsmath,amssymb,amsfonts}
\usepackage{amsthm}

\theoremstyle{plain}
\newtheorem{theorem}{Theorem}[section]

\newtheorem{proposition}{Proposition}[section]
\newtheorem{corollary}{Corollary}[section]

\theoremstyle{definition}
\newtheorem{definition}{Definition}[section]
\newtheorem{example}{Example}[section]
\newtheorem{observation}[theorem]{Observation}

\theoremstyle{remark}

\DeclareMathOperator{\DM}{DM}
\DeclareMathOperator{\Hom}{Hom}
\DeclareMathOperator{\End}{End}
\DeclareMathOperator{\Nat}{Nat}
\DeclareMathOperator{\EM}{EM}
\DeclareMathOperator{\Kl}{Kl}

\title[The Dedekind-MacNeille Completion as a Functor]{The Dedekind-MacNeille Completion as a Functor}
\author{Frank Murphy-Hernandez}
\email{murphy@ciencias.unam.mx}

\begin{document}

\begin{abstract}
We study the Dedekind-MacNeille completion through representable functors on the wide subcategory $\mathbf{Pos}_{\mathrm{cl}}$ of posets and lower cut-continuous (closed) maps. On this category the completion is contravariantly functorial and naturally isomorphic to $\Hom_{\mathbf{Pos}_{\mathrm{cl}}}(-,[1])^{\mathrm{op}}$, via $A\mapsto\chi_A$ and its inverse $f\mapsto f^{-1}(\{0\})$. The naturality of this identification makes the Yoneda lemma applicable: uniqueness of the representing object and the endomorphism monoid of the completion follow formally. We also show that $\mathbf{CLat}_\vee$ is full in $\mathbf{Pos}_{\mathrm{cl}}$, that $\mathbf{CLat}_\wedge$ is equivalent to a full subcategory of $\mathbf{Pos}_{\mathrm{cl}}^{\mathrm{op}}$, and that the covariant adjunction refines to a $2$-adjunction whose idempotent monad has Eilenberg-Moore algebras equivalent to $\mathbf{CLat}_\vee$.
\end{abstract}

\maketitle


\section{Introduction}\label{sec:intro}

It is natural to ask whether every poset can be extended to a lattice. In 1937, H. MacNeille gave a construction embedding any poset in a complete lattice by means of Dedekind cuts \cite{macneille1937partially}. Banaschewski and Bruns described this Dedekind-MacNeille completion as the injective hull in the category of posets with respect to order embeddings \cite{banaschewski1967categorical}. The reflection of complete lattices into posets via lower cut-continuous maps is classical and is due to Bishop \cite{bishop1978universal} and Ern\'e \cite{erne1991dedekind}; we recall it here for convenience. Since $\mathbf{Pos}_{\mathrm{cl}}$ is the classical category of posets and lower cut-continuous maps, its injective-hull character pertains to $\mathbf{Pos}$ with respect to order embeddings. We do not investigate injectivity within $\mathbf{Pos}_{\mathrm{cl}}$ itself; see Corollary~\ref{cor:eta-injective-hull} for a discussion of this distinction.

In analogy with functional analysis, we define the continuous dual of a poset as the subposet of the dual consisting of closed maps, namely monotone maps whose inverse image of a closed set is closed. These are precisely Ern\'e's lower cut-continuous maps.

Our main observation is the natural isomorphism $\DM(X)\cong\Hom_{\mathbf{Pos}_{\mathrm{cl}}}(X,[1])^{\mathrm{op}}$ (Theorem~\ref{thm:macneille-iso}). This yields the functoriality of the completion on the wide subcategory $\mathbf{Pos}_{\mathrm{cl}}$, and, through the Yoneda lemma, formal consequences such as the uniqueness of the representing object, the computation of the endomorphism monoid of the completion, the objectwise idempotence of the double dual, the failure of the canonical comparison to be an isomorphism of presheaves on $\mathbf{Pos}$, and the adjunction between the covariant completion and the inclusion of complete lattices. The adjunction itself is classical; we state and prove it for the reader's convenience. We also record a $2$-categorical refinement of this adjunction, which is not new as a statement about the classical reflection, but which unifies the adjunction, reflexivity, and fullness of $\mathbf{CLat}_\vee$. In \S\ref{subsec:reflection} we characterize the choice of closed maps by proving that $\mathbf{CLat}_\vee$ is a full subcategory of $\mathbf{Pos}_{\mathrm{cl}}$ and that $\mathbf{CLat}_\wedge$ is equivalent to a full subcategory of $\mathbf{Pos}_{\mathrm{cl}}^{\mathrm{op}}$; in \S\ref{sec:monad} we study the induced idempotent monad; and in \S\ref{sec:examples} we collect examples illustrating distributivity, non-injectivity on posets, and cartesian products.

Finally, the use of $[1]=\{0,1\}$ as both the classifying object for MacNeille-closed subsets and the object defining $(-)^\circ$ is informally a restricted ``schizophrenic'' role in Wraith's sense \cite{wraith1993using}; we do not pursue that general framework here.

\section{Preliminaries}\label{sec:prelim}

We denote by $\mathbf{Pos}$ the category of posets and monotone maps; if $X$ is a poset, then $X^{op}$ denotes its opposite poset. A monotone map $f\colon X\longrightarrow Y$ is an embedding if it is injective and $x\leq y$ whenever $f(x)\leq f(y)$. For $A\subseteq X$, we write $A^{\uparrow}$ for the set of upper bounds of $A$ and $A^{\downarrow}$ for the set of lower bounds of $A$; this notation is used uniformly throughout the paper. We say that $A$ is closed if $A^{\uparrow\downarrow}=A$. The operator $A\mapsto A^{\uparrow\downarrow}$ is a closure operator on $\mathbb{P}(X)$: we have $A\subseteq A^{\uparrow\downarrow}$, $A^{\uparrow\downarrow\uparrow\downarrow}=A^{\uparrow\downarrow}$, and $A\subseteq B$ implies $A^{\uparrow\downarrow}\subseteq B^{\uparrow\downarrow}$; its fixed points are precisely the closed subsets. For $x\in X$, the set $\downarrow x=\{y\in X\mid y\leq x\}$ is closed. The Dedekind-MacNeille completion $\DM(X)$ of $X$ is the poset of closed subsets of $X$, ordered by inclusion, and there is a canonical embedding $\epsilon_X\colon X\longrightarrow\DM(X)$ given by $\epsilon_X(x)=\downarrow x$.

With the inclusion order, $\DM(X)$ is a complete lattice: intersections of closed sets are closed, and suprema are obtained by closing the union. The image of $X$ is join-dense and meet-dense in $\DM(X)$, and $\DM(X)$ is the smallest complete lattice with this property; in particular, if $X$ is already a complete lattice, then $\DM(X)\cong X$. For the standard categorical facts about $\DM(X)$ and for the injectivity of complete lattices in $\mathbf{Pos}$, see \cite{banaschewski1967categorical,adamek2004abstract}; for basic order and lattice concepts, see \cite{davey2002introduction}.

We denote by $\mathbf{CLat}_\vee$ the category of complete lattices and suprema-preserving maps, and by $\mathbf{CLat}_\wedge$ the category of complete lattices and infima-preserving maps.

We write $[1]=\{0,1\}$ with its canonical order. We denote by $\mathbf{Pos}_{\mathrm{cl}}$ the category of posets and closed maps: a monotone map $f\colon X\longrightarrow Y$ is closed if $f^{-1}(C)$ is closed in $X$ for every closed $C\subseteq Y$. These are precisely Ern\'e's \emph{lower cut-continuous} maps.

We recall briefly that a presheaf on a category $\mathcal{C}$ is a functor $F\colon\mathcal{C}^{op}\longrightarrow\mathbf{Set}$, and that $F$ is representable if $F\cong\Hom_{\mathcal{C}}(-,A)$ for some object $A\in\mathcal{C}$; equivalently, if it has a universal object. We shall use the equivalence with the Yoneda lemma; for details, see \cite{vistoli2005notes}.

\section{Dual Posets}\label{sec:dual}

The notion of a closed map between posets was recalled in \S\ref{sec:prelim}; it is due to Ern\'e \cite{erne1991dedekind} (see also Bishop \cite{bishop1978universal} for the lattice case). Not every monotone map is closed: for instance, $f\colon\mathbb{R}\to[1]$ with $f(x)=1$ for $x\geq\pi$ and $f(x)=0$ for $x<\pi$ is monotone but not closed, since $f^{-1}(0)=(-\infty,\pi)$ satisfies $(-\infty,\pi)^{\uparrow\downarrow}=(-\infty,\pi]\neq(-\infty,\pi)$.

\begin{definition}\label{def:continuous-dual}
The \emph{continuous dual} of a poset $X$ is $X^*:=\Hom_{\mathbf{Pos}_{\mathrm{cl}}}(X,[1])$, endowed with the pointwise order. For $f,g\in X^*$ we have $f\leq g$ iff $f^{-1}(1)\subseteq g^{-1}(1)$ iff $g^{-1}(0)\subseteq f^{-1}(0)$. The assignment $X\mapsto X^*$ and $u\mapsto u^*:=(-)\circ u$ defines a contravariant functor $(-)^*\colon\mathbf{Pos}_{\mathrm{cl}}^{op}\to\mathbf{Pos}_{\mathrm{cl}}$.
\end{definition}

\begin{definition}\label{def:dm-dual}
The \emph{Dedekind-MacNeille dual} of $X$ is $X^\circ:=(X^*)^{op}$, ordered by $\preceq$. We write $\bigvee^\circ$ and $\bigwedge^\circ$ for suprema and infima in $X^\circ$; thus $f\preceq g$ iff $f^{-1}(0)\subseteq g^{-1}(0)$. The assignment $X\mapsto X^\circ$, $u\mapsto u^\circ:=u^*$, defines a contravariant functor $(-)^\circ\colon\mathbf{Pos}_{\mathrm{cl}}^{\mathrm{op}}\to\mathbf{CLat}_{\wedge}$.
\end{definition}

That $X^*$ is a complete lattice will follow from Theorem~\ref{thm:macneille-iso}; we therefore omit a separate proof. For $\alpha\colon X\to Y$ closed and $\mathcal{A}\subseteq Y^*$, writing $Z(f):=f^{-1}(0)$ one has $Z(\bigvee_{f\in\mathcal{A}}f)=\bigcap_{f\in\mathcal{A}}Z(f)$, and since $\alpha^{-1}$ preserves intersections, $\alpha^*$ preserves suprema; as $X^*$ and $Y^*$ are complete lattices, $\alpha^*$ is closed (Proposition~\ref{prop:closed-iff-sup}), so $(-)^*$ is well defined. Since $X^\circ=(X^*)^{op}$ and taking opposites reverses the direction of the order and exchanges suprema with infima, the same argument shows that $(-)^\circ$ is well defined as a functor into $\mathbf{CLat}_\wedge$.

\begin{definition}\label{def:fx}
For a poset $X$ and $x\in X$, define $f_x\colon X\to[1]$ by $f_x(y)=0$ if $y\leq x$ and $f_x(y)=1$ otherwise.
\end{definition}

Each $f_x$ is closed: it is monotone, and $f_x^{-1}(0)=\downarrow x$ is closed. Moreover, if $f\in X^\circ$ and $f(x)=0$, then $\downarrow x\subseteq f^{-1}(0)$, so $f_x\preceq f$.

\begin{definition}\label{def:eta}
Define $\eta_X\colon X\to X^\circ$ by $\eta_X(x)=f_x$.
\end{definition}

The map $\eta_X$ is monotone, since $x\leq y$ implies $\downarrow x\subseteq\downarrow y$. It is closed: for a closed $\mathcal{A}\subseteq X^\circ$ one has $\eta_X^{-1}(\mathcal{A})=(\bigvee^\circ\mathcal{A})^{-1}(0)$, which is closed because $\bigvee^\circ\mathcal{A}\in X^\circ$. It is injective: if $f_x=f_y$, then $x\in\downarrow x=\downarrow y$ and $y\in\downarrow y=\downarrow x$, so $x=y$. Hence $\eta_X$ is an injective closed map. Under the identification of Theorem~\ref{thm:macneille-iso}, $\eta_X$ corresponds to the canonical embedding $\epsilon_X$ of \S\ref{sec:prelim}. In particular, if $X$ is a complete lattice, then $\eta_X$ is an isomorphism in $\mathbf{Pos}_{\mathrm{cl}}$, since every closed $f\in X^\circ$ satisfies $f=f_{\bigvee f^{-1}(0)}$; consequently $\eta_{X^\circ}\colon X^\circ\to X^{\circ\circ}$ is an isomorphism for every poset $X$.

\section{The MacNeille Isomorphism}\label{sec:macneille-iso}

In this section we prove that the Dedekind-MacNeille dual and the Dedekind-MacNeille completion are naturally isomorphic with respect to closed maps. This identification absorbs the relevant part of \S\ref{sec:dual}, turns the injectivity result into Corollary~\ref{cor:eta-injective-hull}, and streamlines the representability argument of Section~\ref{sec:repr}. For posets $X$ and $Y$ and a closed map $u\colon X\to Y$, we write $\DM(u)\colon\DM(Y)\to\DM(X)$ for the map defined by $\DM(u)(A)=u^{-1}(A)$ for $A\in\DM(Y)$; this is well defined since $u$ is closed.

\begin{theorem}[MacNeille Isomorphism]\label{thm:macneille-iso}
For every poset $X$, the map
\[
\varphi_X\colon X^\circ\longrightarrow\DM(X),\qquad \varphi_X(f):=f^{-1}(\{0\}),
\]
is an isomorphism of complete lattices and an isomorphism in $\mathbf{CLat}_\wedge$. Its inverse is
\[
\psi_X\colon\DM(X)\longrightarrow X^\circ,\qquad \psi_X(A):=\chi_A,
\]
where $\chi_A(x)=0$ if $x\in A$ and $\chi_A(x)=1$ otherwise. Moreover, the family $\{\varphi_X\}_{X\in\mathbf{Pos}_{\mathrm{cl}}}$ is natural: for every closed map $u\colon X\to Y$, the diagram
\[
\begin{array}{ccc}
Y^\circ & \xrightarrow{\ \varphi_Y\ } & \DM(Y)\\[2pt]
{\scriptstyle u^\circ}\big\downarrow & & \big\downarrow{\scriptstyle \DM(u)}\\[2pt]
X^\circ & \xrightarrow{\ \varphi_X\ } & \DM(X)
\end{array}
\]
commutes. Consequently, $\varphi\colon(-)^\circ\Rightarrow\DM$ is a natural isomorphism of functors $\mathbf{Pos}_{\mathrm{cl}}^{\mathrm{op}}\to\mathbf{CLat}_\wedge$.
\end{theorem}

\begin{proof}
\emph{Well-defined maps and mutual inverses.} If $f\in X^\circ$, then $f$ is a closed map $X\to[1]$, so $f^{-1}(\{0\})$ is a closed subset of $X$; hence $\varphi_X(f)\in\DM(X)$. Conversely, let $A\in\DM(X)$ and define $\chi_A\colon X\to[1]$ by $\chi_A(x)=0$ if $x\in A$ and $\chi_A(x)=1$ otherwise. Since $A$ is a lower set, $\chi_A$ is monotone: if $x\leq y$ and $\chi_A(y)=0$, then $y\in A$, hence $x\in A$ and $\chi_A(x)=0$. Since $\chi_A^{-1}(\{0\})=A$ is closed and $\chi_A^{-1}([1])=X$ is closed, $\chi_A$ is a closed map; thus $\psi_X(A)=\chi_A\in X^\circ$. The two assignments are mutually inverse: $\varphi_X(\psi_X(A))=\chi_A^{-1}(\{0\})=A$ for $A\in\DM(X)$, and $\psi_X(\varphi_X(f))=\chi_{f^{-1}(\{0\})}=f$ for $f\in X^\circ$, since $f$ is precisely the characteristic function of its zero set $f^{-1}(\{0\})$.

\emph{Isomorphism of complete lattices.} By Definition~\ref{def:dm-dual}, for $f,g\in X^\circ$ we have $f\preceq g$ if and only if $f^{-1}(0)\subseteq g^{-1}(0)$, which is exactly $\varphi_X(f)\subseteq\varphi_X(g)$. Hence $\varphi_X$ is order-preserving and order-reflecting. Being bijective, order-preserving, and order-reflecting between complete lattices, $\varphi_X$ is an isomorphism of complete lattices and preserves all suprema and infima. In particular, $\varphi_X$ is a morphism in $\mathbf{CLat}_\wedge$, and so is its inverse $\psi_X$.

\emph{Naturality.} Let $u\colon X\to Y$ be a closed map and let $g\in Y^\circ$. Then $u^\circ(g)=g\circ u$, and
\[
\varphi_X(u^\circ(g))
=(g\circ u)^{-1}(0)
=u^{-1}(g^{-1}(0))
=u^{-1}(\varphi_Y(g))
=\DM(u)(\varphi_Y(g)).
\]
Hence $\varphi_X\circ u^\circ=\DM(u)\circ\varphi_Y$, which is precisely the commutativity of the displayed diagram.
\end{proof}

\begin{corollary}\label{cor:varphi-eta-epsilon}
For every poset $X$, with $\eta_X\colon X\to X^\circ$ defined by $\eta_X(x)(y)=0$ if $y\leq x$ and $\eta_X(x)(y)=1$ otherwise, the identity $\varphi_X\circ\eta_X=\epsilon_X$ holds. Consequently, $\eta_X$ is an embedding of posets, and $\eta_X(X)$ is both join-dense and meet-dense in $X^\circ$.
\end{corollary}

\begin{proof}
For $x\in X$, $\eta_X(x)^{-1}(0)=\downarrow x$, so $\varphi_X(\eta_X(x))=\downarrow x=\epsilon_X(x)$. Since $\epsilon_X$ is an embedding and $\varphi_X$ is an isomorphism of posets, $\eta_X$ is an embedding. Since $\epsilon_X(X)$ is join-dense and meet-dense in $\DM(X)$ and $\varphi_X$ is an isomorphism of complete lattices, it maps join-dense and meet-dense subsets to join-dense and meet-dense subsets; hence $\eta_X(X)$ is join-dense and meet-dense in $X^\circ$.
\end{proof}

As a formal consequence of Theorem~\ref{thm:macneille-iso}, the functor $\DM\colon\mathbf{Pos}_{\mathrm{cl}}^{\mathrm{op}}\to\mathbf{CLat}_\wedge$ is naturally isomorphic to $X\mapsto\Hom_{\mathbf{Pos}_{\mathrm{cl}}}(X,[1])^{\mathrm{op}}$, so its underlying set-valued presheaf is representable with representing object $[1]$. We shall use this representability in Section~\ref{sec:repr} without repeating the argument.

\begin{corollary}[Injectivity as a corollary of the MacNeille isomorphism]\label{cor:eta-injective-hull}
Let $X$ be a poset. Under the isomorphism $\varphi_X\colon X^\circ\to\DM(X)$ of Theorem~\ref{thm:macneille-iso}, the embedding $\eta_X\colon X\to X^\circ$ corresponds to the canonical MacNeille embedding $\epsilon_X\colon X\to\DM(X)$, that is, $\varphi_X\circ\eta_X=\epsilon_X$. Consequently, $\epsilon_X$ is an injective hull of $X$ in $\mathbf{Pos}$ with respect to order embeddings, by the classical theorem of Banaschewski and Bruns \cite{banaschewski1967categorical,adamek2004abstract}. Equivalently, $\eta_X$ is an injective hull of $X$ in $\mathbf{Pos}$. This statement concerns $\mathbf{Pos}$ and monotone maps; it does not assert that $\epsilon_X$ or $\eta_X$ is an injective hull in $\mathbf{Pos}_{\mathrm{cl}}$, where the relevant embeddings and injective objects are those of $\mathbf{Pos}_{\mathrm{cl}}$. In particular, the completeness of $\DM(X)$ and of $X^\circ$ makes them injective in $\mathbf{Pos}$, not automatically in $\mathbf{Pos}_{\mathrm{cl}}$.
\end{corollary}

\begin{proof}
The identity $\varphi_X\circ\eta_X=\epsilon_X$ is Corollary~\ref{cor:varphi-eta-epsilon}. Since $\varphi_X$ is an isomorphism of complete lattices, it preserves join-density and meet-density, so $\eta_X(X)$ is join-dense and meet-dense in $X^\circ$. By the classical criterion of Banaschewski and Bruns \cite{banaschewski1967categorical}, an order embedding whose image is both join-dense and meet-dense is essential; hence $\eta_X$ is essential. Since $X^\circ$ is a complete lattice, it is injective in $\mathbf{Pos}$ \cite{banaschewski1967categorical,adamek2004abstract}. Therefore $\eta_X$ is an injective hull in $\mathbf{Pos}$. The same argument applies to $\epsilon_X$ and $\DM(X)$, and the identification $\varphi_X\circ\eta_X=\epsilon_X$ transfers the statement from $\eta_X$ to $\epsilon_X$.
\end{proof}

\section{Representability of the Dedekind-MacNeille Completion}\label{sec:repr}

We now recast the MacNeille isomorphism as a representability statement and derive its formal consequences. The mathematical content is that of Theorem~\ref{thm:macneille-iso}; the Yoneda lemma supplies the universal element, uniqueness, and the endomorphism monoid without a second construction.

\subsection{Representability and Yoneda}\label{subsec:yoneda}

\begin{theorem}[Yoneda formulation of representability]\label{thm:yoneda-repr}
The presheaf $\DM\colon\mathbf{Pos}_{\mathrm{cl}}^{\mathrm{op}}\to\mathcal{S}$ is representable, with representing object $[1]$. The natural isomorphism
\[
\Psi\colon\DM\;\xrightarrow{\ \sim\ }\;\Hom_{\mathbf{Pos}_{\mathrm{cl}}}(-,[1])
\]
is given, for each poset $X$ and each $A\in\DM(X)$, by
\[
\Psi_X(A):=f_A,\qquad
f_A(x)=\begin{cases}0,&x\in A,\\ 1,&x\notin A,\end{cases}
\]
with inverse $\Psi_X^{-1}(f)=f^{-1}(\{0\})$. The pair $([1],\{0\})$ is a universal element of $\DM$: for every $X$ and every $A\in\DM(X)$ there is a unique closed map $f\colon X\to[1]$ with $\DM(f)(\{0\})=A$. Consequently, any other representing object $R$ is canonically isomorphic to $[1]$, by a unique isomorphism carrying the universal element of $R$ to $\{0\}$.
\end{theorem}

\begin{proof}
By Theorem~\ref{thm:macneille-iso}, for every poset $X$ the map $\varphi_X\colon X^\circ\to\DM(X)$, $\varphi_X(f)=f^{-1}(\{0\})$, is a bijection, and the naturality proved therein says precisely that $\varphi\colon\Hom_{\mathbf{Pos}_{\mathrm{cl}}}(-,[1])\Rightarrow\DM$ is a natural isomorphism. Setting $\Psi:=\varphi^{-1}$ gives the stated natural isomorphism. By the Yoneda lemma, the universal element is $\Psi_{[1]}^{-1}(\mathrm{id}_{[1]})=\mathrm{id}_{[1]}^{-1}(\{0\})=\{0\}\in\DM([1])$. The uniqueness statement is the standard full faithfulness of the Yoneda embedding.
\end{proof}

\begin{corollary}[Endomorphisms of the completion]\label{cor:end-DM}
There is a canonical isomorphism of monoids
\[
\End(\DM):=\Nat(\DM,\DM)\;\cong\;
\Hom_{\mathbf{Pos}_{\mathrm{cl}}}([1],[1])\;\cong\;[1],
\]
where the monoid structure on $[1]$ is the meet monoid $([1],\wedge,1)$. Explicitly, $\End(\DM)=\{\mathrm{id},e\}$, where $e_X(A)=X$ for every $A\in\DM(X)$; and $e^2=e$. Likewise, 
\[
\Nat(\Hom_{\mathbf{Pos}_{\mathrm{cl}}}(-,[1]),\DM)\cong\DM([1])=\{\{0\},[1]\}
\]
, consisting of the isomorphism $\varphi$ and the constant-at-maximum transformation.
\end{corollary}

\begin{proof}
By full faithfulness of the Yoneda embedding, $\Nat(\Hom(-,[1]),\Hom(-,[1]))\cong\Hom_{\mathbf{Pos}_{\mathrm{cl}}}([1],[1])$. The closed maps $[1]\to[1]$ are $\mathrm{id}_{[1]}$ and the constant map $c_0$ with value $0$; thus the monoid is $[1]$ under meet. Under Yoneda, $\mathrm{id}_{[1]}$ corresponds to the identity of $\DM$, and $c_0$ corresponds to $e_X(A)=\Psi_X^{-1}(c_0\circ\Psi_X(A))=X$, since the zero set of $c_0\circ f_A$ is $X$. The last assertion is the Yoneda bijection $\Nat(\Hom(-,[1]),\DM)\cong\DM([1])$.
\end{proof}

\begin{observation}[Failure of the canonical comparison on $\mathbf{Pos}$]\label{obs:not-repr-pos}
The formula $A\mapsto\chi_A$ defines an injective map of sets $\Phi_X\colon\DM(X)\to\Hom_{\mathbf{Pos}}(X,[1])^{\mathrm{op}}$ for every poset $X$, because $A$ is determined by $\chi_A$. It is not surjective: the monotone map $\mathbb{R}\to[1]$ with jump at $\pi$ has zero set $(-\infty,\pi)$, which is not closed, so it is not in the image of $\Phi_{\mathbb{R}}$. Thus the canonical comparison does not yield a natural isomorphism of presheaves on $\mathbf{Pos}$; Theorem~\ref{thm:yoneda-repr} is genuinely a statement about $\mathbf{Pos}_{\mathrm{cl}}$.
\end{observation}

\subsection{Reflection and adjunction}\label{subsec:reflection}

\begin{proposition}\label{prop:closed-iff-sup}
Let $L$ and $M$ be complete lattices and let $f\colon L\to M$ be monotone. Then $f$ is closed if and only if $f$ preserves arbitrary suprema.
\end{proposition}

\begin{proof}
Recall that every closed subset of a complete lattice is a principal ideal: if $B$ is closed and nonempty, then $B=\downarrow(\bigvee B)$; and $\varnothing$ is not closed because $\varnothing^{ul}=\{\bot\}$.

If $f$ preserves suprema, then for $b\in M$ the set $f^{-1}(\downarrow b)$ is a lower set closed under suprema, hence equals $\downarrow a$ with $a=\bigvee f^{-1}(\downarrow b)$. Thus preimages of principal ideals are principal ideals, so $f$ is closed.

Conversely, suppose $f$ is closed. Let $S\subseteq L$, $a=\bigvee S$, and $b=\bigvee f(S)$. Monotonicity gives $b\leq f(a)$. The principal ideal $\downarrow b$ is closed in $M$, so $f^{-1}(\downarrow b)$ is closed in $L$, hence equals $\downarrow a'$ for some $a'$. Since $f(S)\subseteq\downarrow b$, we have $S\subseteq\downarrow a'$, so $a\leq a'$, and $f(a)\leq f(a')\leq b$. Hence $f(a)=b$, and $f$ preserves suprema.
\end{proof}

\begin{corollary}\label{cor:CLat-vee-full}
$\mathbf{CLat}_\vee$ is a full subcategory of $\mathbf{Pos}_{\mathrm{cl}}$: for $L,M\in\mathbf{CLat}_\vee$,
 \[\Hom_{\mathbf{CLat}_\vee}(L,M)=\Hom_{\mathbf{Pos}_{\mathrm{cl}}}(L,M).\]
\end{corollary}

\begin{corollary}\label{cor:CLat-wedge-full}
The adjoint correspondence on complete lattices yields an isomorphism $\mathbf{CLat}_\wedge^{\mathrm{op}}\cong\mathbf{CLat}_\vee$. Consequently, $\mathbf{CLat}_\wedge$ is equivalent to a full subcategory of $\mathbf{Pos}_{\mathrm{cl}}^{\mathrm{op}}$.
\end{corollary}

\begin{proof}
Each infima-preserving map has a unique left adjoint, which is suprema-preserving, and this correspondence is natural; combined with Corollary~\ref{cor:CLat-vee-full}, this gives the claim.
\end{proof}

\begin{theorem}[Adjunction]\label{thm:adjunction}
Let $\DM_{\mathrm{cov}}\colon\mathbf{Pos}_{\mathrm{cl}}\to\mathbf{CLat}_\vee$ be defined on objects by $\DM_{\mathrm{cov}}(X)=\DM(X)$ and on morphisms by $\DM_{\mathrm{cov}}(u)(A):=u(A)^{ul}$ for $u\colon X\to Y$ closed and $A\in\DM(X)$. Let $\iota\colon\mathbf{CLat}_\vee\hookrightarrow\mathbf{Pos}_{\mathrm{cl}}$ be the inclusion. Then $\DM_{\mathrm{cov}}\dashv\iota$. The unit is $\eta_X=\epsilon_X\colon X\to\DM(X)$, the counit is $\varepsilon_L(A)=\bigvee_L A$, and for $X\in\mathbf{Pos}_{\mathrm{cl}}$, $L\in\mathbf{CLat}_\vee$, the map
\[
\Phi_{X,L}\colon\Hom_{\mathbf{CLat}_\vee}(\DM(X),L)\longrightarrow\Hom_{\mathbf{Pos}_{\mathrm{cl}}}(X,L),\qquad \Phi_{X,L}(g):=g\circ\epsilon_X,
\]
is a natural bijection with inverse
\[
\Psi_{X,L}(h)(A):=\bigvee_{x\in A}h(x),\qquad A\in\DM(X).
\]
This adjunction is classical (Bishop \cite{bishop1978universal}, Ern\'e \cite[\S 2]{erne1991dedekind}); we include the proof to fix conventions.
\end{theorem}

\begin{proof}
\emph{Well defined.} If $g\in\mathbf{CLat}_\vee$, then $g$ is closed by Proposition~\ref{prop:closed-iff-sup}; the embedding $\epsilon_X$ is closed, since the preimage of a principal ideal $\downarrow A$ is $A$. Hence $g\circ\epsilon_X$ is closed. Conversely, if $h\colon X\to L$ is closed, then $\Psi(h)(A)=\bigvee_{x\in A}h(x)$ is the unique suprema-preserving extension of $h$ along the join-dense embedding $\epsilon_X$, as shown below.

\emph{Mutual inverses.} For $h\colon X\to L$ closed, $\Psi(h)(\epsilon_X(x))=\bigvee_{y\in x_\leq}h(y)=h(x)$. Conversely, if $g\in\mathbf{CLat}_\vee$, then by join-density of $\epsilon_X(X)$ and sup-preservation of $g$,
\[
g(A)=g\Bigl(\bigvee_{x\in A}\epsilon_X(x)\Bigr)=\bigvee_{x\in A}g(\epsilon_X(x))=\Psi(g\circ\epsilon_X)(A).
\]

\emph{$\Psi(h)$ preserves suprema.} Let $\mathcal{A}\subseteq\DM(X)$ and $A=\bigvee^{\DM(X)}\mathcal{A}=(\bigcup\mathcal{A})^{ul}$. It suffices to show, for every $l\in L$,
\[
\bigvee_{x\in A}h(x)\leq l \quad\Longleftrightarrow\quad \bigvee_{B\in\mathcal{A}}\bigvee_{x\in B}h(x)\leq l.
\]
The forward implication is immediate from $\bigcup\mathcal{A}\subseteq A$. Conversely, if $h(x)\leq l$ for all $x\in\bigcup\mathcal{A}$, then $\bigcup\mathcal{A}\subseteq h^{-1}(\downarrow l)$; since $h$ is closed, $h^{-1}(\downarrow l)$ is closed in $X$, hence contains $(\bigcup\mathcal{A})^{ul}=A$, and therefore $\bigvee_{x\in A}h(x)\leq l$. The empty family is covered by the same argument: $A=\varnothing^{ul}$ is the least closed subset, and the hypothesis is vacuous.

\emph{Naturality.} For $u\colon X'\to X$ closed and $g\colon\DM(X)\to L$,
\[
\Phi_{X',L}(\DM_{\mathrm{cov}}(u)^*(g))=(g\circ\DM_{\mathrm{cov}}(u))\circ\epsilon_{X'}=g\circ(\epsilon_X\circ u)=(g\circ\epsilon_X)\circ u=u^*(\Phi_{X,L}(g)),
\]
using $\DM_{\mathrm{cov}}(u)\circ\epsilon_{X'}=\epsilon_X\circ u$, i.e.\ $u(x'_\leq)^{ul}=u(x')_\leq$. Naturality in $L$ follows from sup-preservation of morphisms in $\mathbf{CLat}_\vee$. The triangle identities follow from the mutual-inverse calculation above, with $\varepsilon_L\circ\eta_L=\mathrm{id}_L$ and $\eta_L\circ\varepsilon_L=\mathrm{id}_{\DM(L)}$, where $\eta_L := \epsilon_L$ is the unit.
\end{proof}

\begin{corollary}[Reflexivity]\label{cor:reflexive}
$\mathbf{CLat}_\vee$ is a full reflective subcategory of $\mathbf{Pos}_{\mathrm{cl}}$, with reflector $\DM_{\mathrm{cov}}$ and counit $\varepsilon_L(A)=\bigvee_L A$. The counit is an isomorphism with inverse $x\mapsto x_\leq$, the fixed points of $\DM_{\mathrm{cov}}$ are exactly the complete lattices, and the adjunction restricts on $\mathbf{CLat}_\vee$ to the identity.
\end{corollary}

\begin{proof}
The adjunction is Theorem~\ref{thm:adjunction} and fullness is Corollary~\ref{cor:CLat-vee-full}. The counit is an isomorphism because $L$ is complete, with the stated inverse by the triangle identity. The fixed points are the complete lattices by the standard characterization of the MacNeille completion; equivalently, the unit $\epsilon_X$ is an isomorphism exactly when $X$ is complete. On $\mathbf{CLat}_\vee$ the counit identifies $\DM(L)$ with $L$, so the restriction is the identity.
\end{proof}

\begin{proposition}[Maximality of $\mathbf{Pos}_{\mathrm{cl}}$]\label{obs:pos-cl-maximal}
Let $\mathcal{C}$ be a wide subcategory of $\mathbf{Pos}$ such that for every poset $X$ the canonical map $\DM(X)\to\Hom_{\mathcal{C}}(X,[1])$, $A\mapsto\chi_A$, is a bijection. Then $\mathcal{C}\subseteq\mathbf{Pos}_{\mathrm{cl}}$. Since $\mathbf{Pos}_{\mathrm{cl}}$ itself has this property by Theorem~\ref{thm:yoneda-repr}, it is the largest wide subcategory with this property. This is a reformulation of the known maximality of lower cut-continuous maps (Ern\'e \cite{erne1991dedekind}); the short argument is included for completeness.
\end{proposition}

\begin{proof}
Let $u\colon X\to Y$ be a morphism of $\mathcal{C}$ and let $A\subseteq Y$ be closed. Then $\chi_A\in\Hom_{\mathcal{C}}(Y,[1])$, hence $\chi_A\circ u\in\Hom_{\mathcal{C}}(X,[1])$. By hypothesis there is a unique closed $B\in\DM(X)$ with $\chi_A\circ u=\chi_B$. Taking zero sets gives $u^{-1}(A)=B$, which is closed. Thus $u$ is closed, so $\mathcal{C}\subseteq\mathbf{Pos}_{\mathrm{cl}}$.
\end{proof}


\section{Examples and Computations}\label{sec:examples}

We retain only the examples that illustrate distinct mathematical phenomena: the appearance of new points in the completion, the failure of distributivity, and the failure of the completion to determine the original poset.

\begin{example}[The rationals]\label{ex:Q}
The closed subsets of $\mathbb{Q}$ are $\varnothing$, $\mathbb{Q}$, the principal ideals $\downarrow q$ for $q\in\mathbb{Q}$, and the lower sets determined by irrational Dedekind cuts. Hence $\DM(\mathbb{Q})$ is order-isomorphic to the extended real line $\overline{\mathbb{R}}=\mathbb{R}\cup\{-\infty,+\infty\}$, with $\varnothing\leftrightarrow-\infty$, $\downarrow q\leftrightarrow q$, $\mathbb{Q}\leftrightarrow+\infty$, and each irrational cut $\alpha$ corresponding to $\alpha\in\mathbb{R}\setminus\mathbb{Q}$. This exhibits the new points produced by the completion and recovers Dedekind's construction of the real numbers from the rationals. By Theorem~\ref{thm:macneille-iso}, $\mathbb{Q}^\circ\cong\DM(\mathbb{Q})\cong\overline{\mathbb{R}}$.
\end{example}

The same phenomenon is visible in simpler chains: $\DM(\mathbb{N})=\mathbb{N}\cup\{+\infty\}$, $\DM(\mathbb{Z})=\mathbb{Z}\cup\{-\infty,+\infty\}$, while for a finite chain $[n]$ the completion is $[n]$ itself because $[n]$ is already complete. In general, for a chain $C$, the poset $\DM(C)$ is the Dedekind completion of $C$, obtained by adjoining the missing extremal endpoints and one new point for each non-principal MacNeille cut.

\begin{example}[The three-element antichain]\label{ex:antichain3}
Let $X=\{a,b,c\}$ be an antichain. Its closed subsets are $\varnothing$, the three singletons $\{a\}$, $\{b\}$, $\{c\}$, and $X$ itself. Hence $\DM(X)=M_3$, the lattice with minimum, three atoms, and maximum. Since $M_3$ is the smallest non-distributive lattice, the MacNeille completion of a poset need not be distributive; in particular, the adjective ``complete lattice'' in the codomain of $\DM$ cannot be sharpened to ``distributive lattice.''
\end{example}

\begin{example}[Failure of injectivity]\label{ex:V-wedge}
Let \(V=\{a,b,c\}\) with \(a<c\) and \(b<c\). The closed subsets of \(V\) are \(\varnothing\), \(\{a\}\), \(\{b\}\), and \(V\), so \(\DM(V)\) is the diamond \([1]\times[1]\). Now let \(W=\{a,b,c,d\}\) be the four-element Boolean algebra, with \(a,b\) minimal and \(c,d\) maximal. Since \(W\) is already a complete lattice, \(\DM(W)\cong W\), which is also the diamond \([1]\times[1]\). Thus \(\DM(V)\cong\DM(W)\), although \(V\) and \(W\) are not isomorphic (they have different cardinalities). Therefore \(\DM\) is not injective on finite posets: the completion does not determine the original poset.
\end{example}

The contrast with the complete case is instructive: on complete lattices, $\DM$ restricts to the identity up to canonical isomorphism, so no information is lost there. The loss of information in Example~\ref{ex:V-wedge} is therefore concentrated in the non-complete part of the poset.

\section{The MacNeille Monad}\label{sec:monad}

The reflection of Theorem~\ref{thm:adjunction} upgrades canonically to a $2$-categorical statement, yielding an idempotent monad whose algebras are the complete lattices. We state only the structural facts and refer to the standard theory of $2$-categories, $2$-monads, and idempotent monads \cite{lack2010two,blackwell1989two,borceux1994handbook}. The bireflection of preorders along lower cut-continuous maps, in the enriched setting, already appears in Garner \cite{garner2014topological}; the $2$-categorical packaging below records the same reflection in the notation of the present paper, and its contribution here is organizational rather than foundational.

\subsection{The $2$-categorical structure and the $2$-adjunction}

\begin{proposition}\label{prop:macneille-2-structure}
The following hold.
\begin{enumerate}
\item $\mathbf{Pos}_{\mathrm{cl}}$, $\mathbf{CLat}_\vee$, and $\mathbf{CLat}_\wedge$ are $2$-categories, with $2$-cells $u\Rightarrow v$ given by the pointwise order $u\leq v$.
\item The functor $\DM_{\mathrm{cov}}\colon\mathbf{Pos}_{\mathrm{cl}}\to\mathbf{CLat}_\vee$ is a $2$-functor.
\item The adjunction $\DM_{\mathrm{cov}}\dashv\iota$ of Theorem~\ref{thm:adjunction} is a $2$-adjunction: the bijection $\Phi_{X,L}$ is an isomorphism of Hom-posets.
\end{enumerate}
\end{proposition}

\begin{proof}
(1) The $2$-category axioms follow from the pointwise order, monotonicity of composition, and associativity of composition of closed maps. (2) If $u\leq v$ pointwise, then for $A\in\DM(X)$ we have $u(A)\subseteq\downarrow v(A)$, hence $u(A)^{ul}\subseteq(\downarrow v(A))^{ul}=v(A)^{ul}$, so $\DM_{\mathrm{cov}}(u)\leq\DM_{\mathrm{cov}}(v)$ pointwise. (3) The maps $\Phi_{X,L}(g)=g\circ\epsilon_X$ and $\Psi_{X,L}(h)(A)=\bigvee_{x\in A}h(x)$ are monotone with respect to the pointwise orders: $g\leq g'$ implies $g\circ\epsilon_X\leq g'\circ\epsilon_X$, and $h\leq h'$ implies $\bigvee_{x\in A}h(x)\leq\bigvee_{x\in A}h'(x)$. Since they are mutually inverse bijections (Theorem~\ref{thm:adjunction}), they are isomorphisms of Hom-posets.
\end{proof}

\subsection{The MacNeille monad}

\begin{definition}\label{def:macneille-monad}
Let $T:=\iota\circ\DM_{\mathrm{cov}}\colon\mathbf{Pos}_{\mathrm{cl}}\to\mathbf{Pos}_{\mathrm{cl}}$. Explicitly, $T(X)=\DM(X)$ on objects, $T(u)(A)=u(A)^{ul}$ on closed maps $u\colon X\to Y$ and $A\in\DM(X)$, the unit is $\eta_X=\epsilon_X\colon X\to\DM(X)$, and the multiplication is $\mu_X\colon\DM(\DM(X))\to\DM(X)$, $\mu_X(\mathcal{A})=\bigvee_{\DM(X)}\mathcal{A}$.
\end{definition}

The data $(T,\eta,\mu)$ form a monad on $\mathbf{Pos}_{\mathrm{cl}}$ by the standard theory of $2$-adjunctions; the unit laws and associativity are the triangle identities verified in the proof of Theorem~\ref{thm:adjunction}. The monad is idempotent: $T^2\cong T$, since $T(X)=\DM(X)$ is a complete lattice and $\mu_X$ is the counit at $T(X)$, an isomorphism by Corollary~\ref{cor:reflexive}.

\begin{corollary}[Eilenberg-Moore algebras]\label{thm:eilenberg-moore-clat}
The comparison functor $\mathbf{CLat}_\vee\to\EM(T)$ is an equivalence of categories, and the same conclusion holds for the Kleisli category: $\Kl(T)\simeq\mathbf{CLat}_\vee\simeq\EM(T)$. The $T$-algebra structure on each $L\in\mathbf{CLat}_\vee$ is unique, given by the counit $\varepsilon_L(A)=\bigvee_L A$.
\end{corollary}

\begin{proof}
By Corollary~\ref{cor:reflexive}, $\mathbf{CLat}_\vee$ is a full reflective subcategory of $\mathbf{Pos}_{\mathrm{cl}}$. The standard theorem on reflective subcategories \cite[\S IV.3]{borceux1994handbook} yields that the induced monad $T$ is idempotent, the comparison functor $\mathbf{CLat}_\vee\to\EM(T)$ is an equivalence, and the algebra structure on each $L\in\mathbf{CLat}_\vee$ is unique. Idempotent monads have equivalent Kleisli and Eilenberg-Moore categories, whence $\Kl(T)\simeq\mathbf{CLat}_\vee\simeq\EM(T)$.
\end{proof}

The idempotent monad $T$ packages the adjunction, reflexivity, and fullness of $\mathbf{CLat}_\vee$ as three facets of one structure. The bireflection of preorders along lower cut-continuous maps is classical \cite{bishop1978universal}; in the enriched setting, the corresponding bireflection is established by Garner \cite{garner2014topological}. The contribution of the present section is therefore the explicit identification of the $2$-cells and Hom-posets in $\mathbf{Pos}_{\mathrm{cl}}$ and the matching of the idempotent monad with the covariant completion of Theorem~\ref{thm:adjunction}; no claim of originality is made for the underlying bireflection.

The covariant and contravariant completions are related by a Galois connection. For $u\colon X\to Y$ closed, the maps $u_!(-):=u(-)^{ul}\colon\DM(X)\to\DM(Y)$ and $u^*(-):=u^{-1}(-)\colon\DM(Y)\to\DM(X)$ satisfy $u_!\dashv u^*$ with respect to the inclusion orders; thus $u_!$ is a morphism in $\mathbf{CLat}_\vee$ and $u^*$ is a morphism in $\mathbf{CLat}_\wedge$. The functors $\DM_{\mathrm{cov}}$ and $\DM$ are recovered by $u\mapsto u_!$ and $u\mapsto u^*$ respectively. This is the precise sense in which the two completions are dual; it does not amount to a $2$-isomorphism between $\DM$ and $\DM_{\mathrm{cov}}$, since they differ in variance and in target categories.


\bibliographystyle{amsplain}
\bibliography{biblio}

\end{document}